\documentclass[12pt]{amsart}
\usepackage[utf8]{inputenc}
\usepackage{amsmath, amsfonts, amssymb, amsthm} 
\usepackage{amsfonts, graphics,epsfig}
\usepackage{amsmath, amsfonts, amssymb, amsthm, amscd}
\usepackage{hyperref}
\usepackage{graphicx}
\usepackage{multicol,multirow}
\usepackage{amsmath,amssymb,amsfonts}
\usepackage{mathrsfs}
\usepackage{rotating}
\usepackage{appendix}
\usepackage[numbers]{natbib}
\usepackage{amsthm}
\usepackage{lipsum}

\newtheorem{theorem}{Theorem}[section]
\newtheorem{lemma}[theorem]{Lemma}
\newtheorem{corollary}[theorem]{Corollary}
\newtheorem{proposition}[theorem]{Proposition}
\newtheorem{remark}[theorem]{Remark}
\newtheorem{example}[theorem]{Example}

\numberwithin{equation}{section}

\newcommand{\SL}{\ensuremath{\operatorname{SL}}}

\newcommand{\SO}{\operatorname{SO}}

\newtheorem{acknowledgement}{Acknowledgement}

\newcommand{\R}{\ensuremath{\mathbb{R}}}

\def\Ga{{\Gamma}}
\def\Gr{{\rm Gr}}

\def\mP{{\mathbb P}}

\def\dim{\text {dim}}

\date{}

\def\build#1_#2^#3{\mathrel{\mathop{\kern 0pt#1}\limits_{#2}^{#3}}}

\let\ap=\alpha

\def\smallsquare{\vbox{\hrule\hbox{\vrule height 1 ex\kern 1
ex\vrule}\hrule}}

\def\limsup{{\rm limsup}}

\def\cO{{\mathcal O}}

\title{Invariant random subgroups on certain orbits}
\author{Manoj Choudhuri} 

\address{Department of Basic Sciences, Institute of Infrastructure, Technology, Research And Management, Near Khokhra Circle, Maninagar (East), Ahmedabad 380026, Gujarat, India}
\email{manojchoudhuri@iitram.ac.in}
\author{C. R. E. Raja}
\address{Theoretical Statistics and Mathematics Unit, Indian Statistical Institute, 8th Mile, Mysore Road, Bangalore 560059, Karnataka, India}
 \email{creraja@isibang.ac.in}

\begin{document}
\maketitle
\begin{abstract}
 Let $G$ be a connected Lie group and $\text{Sub}_G$ be the space of closed subgroups of $G$ equipped with the Chabauty topology. In this article, we investigate the existence of invariant random subgroups of $G$ supported on various orbits of the conjugation action of $G$ on $\text{Sub}_G$.
 \end{abstract}
\keywords{Connected Lie groups, algebraic groups, spaces of closed subgroups, Chabauty topology, invariant random subgroups}

\subjclass[2020]{Primary: $28$D$15$, $22$E$15$, Secondary: $37$A$15$, $37$B$02$}
\tableofcontents
\section{\bf Introduction}
  Let $G$ be a connected (real) Lie group and $\text{Sub}_G$ be the space of (closed) subgroups of $G$ equipped with the Chabauty topology introduced by Claude Chabauty in \cite{CC}. The Chabauty topology on $\text{Sub}_G$ is generated by the subsets $\{ H \in {\rm Sub}_G \mid H\cap K = \emptyset \}$, $\{ H\in {\rm Sub}_G \mid H \cap U \not = \emptyset \}$, 
where $K \subset G$ is compact and $U \subset G$ is open.  With respect to the
Chabauty topology, $\text{Sub}_G$ is a compact metrizable space (cf. Lemma E.1.1 of \cite{BaP}).  Moreover, a sequence $H_n\in\text{Sub}_G$ converges
to $H\in \text{Sub}_G$ if and only if for any $h\in H$ there is a sequence $h_n\in H_n$ such that $h_n$ converges to $h$, and for
any sequence $x_{k_n}\in H_{k_n}$, with $k_{n+1}>k_n$ and $x_{k_n}\to x$, we have $x\in H$ (cf. Lemma E.1.2 of \cite{BaP}). 

The group $G$ acts on $\text{Sub}_G$ by conjugation: $(g, H) \mapsto gHg^{-1}$.  
An invariant random subgroup (hereafter IRS) is a probability measure $\mu$ on $\text{Sub}_G$ that is invariant under the conjugation action of $G$ on $\text{Sub}_G$.  Dirac measure supported on a closed normal subgroup is an IRS.  Also, if $\Gamma$ is a lattice in $G$, i.e., $\Gamma$ is a discrete subgroup of $G$ such that $G/\Gamma$ carries a $G$-invariant probability measure $\mu$, then the push forward of $\mu$ under the map $g$ to $g\Gamma g^{-1}$ defines an IRS on $G$.  Thus, invariant random subgroups can be thought of as a generalization of both normal subgroups and lattices.  The study of invariant random subgroups and ergodic theory are deeply intertwined as an IRS can be realized as stabilizers of a measure-preserving action of the group (see \cite{A7s0} and  \cite{AGV} for more details). 

We say that a Borel subset $\cO$ of $\text{Sub}_G$ supports an IRS $\mu$ if $\mu ({\cO}) =1$:  a Borel subset ${\cO}$ 
supporting an IRS $\mu$ need not be the support of $\mu$ but 
$\overline {\cO}$ contains the support of $\mu$.

In recent times, IRS has attracted many studies (cf.  \cite{A7s0}, \cite{A7s}, \cite{AGV}, \cite{BT}, \cite{GL} and the references therein). In this article, we are interested in the existence of invariant random subgroups supported on certain orbits of the conjugation action of $G$ on $\text{Sub}_G$.  We first consider the orbit of semisimple Levi subgroups.  Let $G$ be a connected Lie group with a closed semisimple Levi subgroup $L$; it may be noted that semisimple Levi subgroups of connected Lie groups need not be closed, but, Levi subgroups of linear groups are closed.  Let $L_G$ be the set of all semisimple Levi subgroups of $G$.  Then $L_G = \{xLx^{-1} \mid x \in N \}$ 
where $N$ is the nilradical of $G$ (cf. Theorem 3.18.3 of \cite{Va}).   
An IRS $\mu$ is called a Levi IRS if $\mu (L_G)=1$.

\begin{theorem}\label{irsl}
Let $G$ be a connected Lie group with nilradical $N$, and $L$ be a semisimple Levi subgroup of $G$.  Then the following are equivalent: 
\begin{enumerate}
   
 \item[(1)] there is an IRS supported on the set of all semisimple Levi subgroups;

\item[(2)] there is a $N$-invariant probability measure supported on semisimple Levi subgroups;

\item[(3)] the semisimple Levi subgroup $L$ is normal in $G$.
\end{enumerate}

In particular, the set of semisimple Levi subgroups $L_G$ is closed in $\text{Sub}_G$ if and only if $L$ is normal. 
\end{theorem}
\begin{corollary}
Any Levi IRS of a connected Lie group $G$ is a Dirac measure supported on a Levi subgroup of $G$.
\end{corollary}

Thus, any connected Lie group $G$ admitting a Levi IRS is of the form $G=L\times S$ where $L$ is semisimple and $S$ is solvable. In particular, for $n \geq 2$, $\SL(n,\R)\times \R^n$ admits a Levi IRS, whereas $\SL(n,\R)\ltimes \R^n$ does not admit a Levi IRS.

As in the case of Levi subgroups, the set of all maximal compact subgroups of a connected Lie group is also the orbit of a maximal compact subgroup under the conjugation action. We next provide a necessary and sufficient condition for the set of maximal compact subgroups to support an IRS.

\begin{theorem}\label{irsK}
Let $G$ be a connected Lie group and $K$ be a maximal compact subgroup of $G$. Also, let $K_G$ be the set of all maximal compact subgroups of $G$. Then,  
there is an IRS supported on $K_G$ if and only if $K$ is normal in $G$.
\end{theorem}

We now assume that $G$ is a real algebraic group, that is, $G$ is the group of $\R$-points of an algebraic group defined over $\R$. Let $B$ be a Borel subgroup of $G$, i.e., $B$ is a maximal connected solvable subgroup of $G$. Also, let $B_G$ be the set of all Borel subgroups of $G$. Then $B_G$ is the orbit of $B$ under the conjugation action. Suppose $H$ is a locally compact group which acts continuously on a locally convex topological vector space $V$. Then, every $H$-invariant, compact, convex subset of $V$ is called a compact, convex $H$-space. We call $H$ to be amenable if $H$ has a fixed point in every compact, convex $H$-space. Abelian, compact and solvable Lie groups are examples of amenable groups. The readers are referred to Chapter $12$ of \cite{MDW} for more details regarding amenable groups.  The next theorem is about the existence of an IRS supported on $B_G$.  

\begin{theorem}\label{irsB}
Let $G$ be a connected real algebraic group and $B_G$ be the set of all Borel subgroups of $G$. Then $B_G$ is closed in $\text{Sub}_G$, and there is an IRS supported on $B_G$ if and only if $G$ is amenable.
\end{theorem}

\begin{example}
Let $G= SO(3, \R) \ltimes \R ^3$.  We identify $O(2, \R)$ as a subgroup of $SO(3, \R)$ by the map 
$\alpha \mapsto \begin{pmatrix} \alpha & 0 \cr 0 & \det (\ap ) \end{pmatrix}$.  Then $B= SO(2, \R) \ltimes R^3$ is a Borel subgroup of $G$ and the normalizer of $B$ is $O(2, \R)\ltimes \R^3$.  Thus, any $SO(3, \R)$-invariant measure on the quotient space $SO(3, \R)/O(2, \R)$ induces an IRS on the set of Borel subgroups $B_G$ of $G$.  
\end{example}

We next consider the orbit of a maximal diagonalizable subgroup $D$. Here, by diagonalizable subgroup, we mean diagonalizable over $\R$, that is, an $\R$-split torus. The set of all maximal diagonalizable subgroups then consists of conjugates of $D$, we denote this set by $D_G$.

\begin{theorem}\label{irsD}
Let $G$ be a connected real algebraic group and $D_G$ be the set of all maximal diagnolizable subgroups of $G$.  Then the following are equivalent: 

\begin{enumerate}
    \item there is an IRS supported on $D_G$;
\item $D$ is central in $G$;
\item $D$ is normal in $G$.
    
\end{enumerate}
Moreover, if (1) holds, then any Levi subgroup of $G$ is compact and the solvable radical is a compact extension of its nilradical.
\end{theorem}

\section{\bf Preliminary}
In this section, we now collect a few results that connect the dynamics of $G$ on $\text{Sub}_G$ with projective linear actions on Grassmannians.  Given a finite-dimensional vector space $V$, let $\Gr _k(V)$ be the Grassnannian manifold of $k$-dimensional subspaces of $V$.  

We include a proof of the following lemma because we are unable to locate it in the following form that compares the Grassmannian topology and the Chabauty topology of $\Gr _k(V)$ seen as a subset of $\text{Sub}_V$. Recall that the metric on $\Gr _k(V)$ is given by  
\[
d(W,W')=||P_W-P_{W'}||
\]  
for $W,W'\in \Gr _k(V)$, with $P_W$ and $P_W'$ being the orthogonal projections onto the subspaces $W$ and $W'$ respectively. One may look at \cite{GL} for the measurability of the map.  

\begin{lemma}\label{irsl2}
Let $V$ be a finite-dimensional vector space and $\Gr_k (V)$ be the Grassmannian manifold of $k$-dimensional subspaces of $V$.   
Let $\mathfrak{i} \colon \Gr _k(V) \to {\rm Sub} _V$ be the canonical inclusion.   
Then $\Gr_k (V)$ is homeomorphic to its image.  In particular, the set of all subspaces of dimension $k$ is a closed set in the Chaubauty topology.    
\end{lemma}

\begin{proof}
Suppose $W_n\to W$ in $\Gr _k(V)$.   This implies that $||P_{W_n}-P_W||\rightarrow 0$ as $n\rightarrow\infty$.
Now let $w\in W$, we want to show that there exist $w_n\in W_n$ such that $w_n\rightarrow w$. Let $w_n=P_{W_n}(w)$. Then
clearly $w_n\in W_n$ and $$||w_n-w||=||P_{W_n}(w)-P_W(w)||\leq ||P_{W_n}-P_W||||w||\rightarrow 0$$ as $n\rightarrow \infty$. On
the other hand assume that there exist $w_n\in W_n$ such that $w_n\rightarrow w$, we want to show that $w\in W$. It is enough to
show that $P_W(w)=w$. Since $w_n\rightarrow w$ and $P_{W_n}\rightarrow P_W$, it follows that $P_{W_n}(w_n)\rightarrow P_W(w)$. But
as $P_{W_n}(w_n)=w_n$, and $w_n\rightarrow w$, it follows that $P_W(w)=w$.  Thus, $W_n \to W$ in Sub$_V$.  This proves that the 
map $\mathfrak{i}$ is continuous.  Since $\Gr _k(V)$ is compact, $\Gr _k(V)$ is homeomorphic to $\mathfrak{i} (\Gr _k(V))$.   
\end{proof}

It is known that if a sequence of closed subgroups $H_n$ of a Lie group $G$ converges to a closed subgroup $H$, then 
dim$(H) \geq \limsup ~ \dim (H_n)$ (cf. Theorem 1 of \cite{Wa}).  However, in the case of closed subgroups of the same dimension, we have the following.

\begin{proposition}\label{cs}
Let $G$ be a connected Lie group and $H_n \to H$ in $\text{Sub}_G$.  Let $\mathcal{ G}$, $\mathcal{H} _n$ and $\mathcal{H}$ be Lie algebras of $G$, $H_n$ and $H$ respectively.    
If $\dim (H_n) = \dim (H)$ for all $n$, then $\mathcal{H} _n \to \mathcal{H}$ in Sub$_\mathcal{G}$. 
\end{proposition}

\begin{remark}

The map sending a closed subgroup to its Lie algebra is known to be a measurable map (cf. proof of Lemma 6 in \cite{Zi}) but here we prove continuity under some restrictions 
though measurability is sufficient for our purpose.  

\end{remark}

\begin{proof}
Let $V$ be a limit of $\mathcal{H}_n$ in Sub$_\mathcal{G}$.  By passing to a subsequence, we may assume that $\mathcal{H}_n \to V$ in 
Sub$_\mathcal{G}$. Then $V$ is a vector subspace, and it is easy to see that $V$ is a Lie subalgebra of $\mathcal{G}$ as well. For $v\in V$, there are $v_n \in \mathcal{H} _n$ such that 
$v_n \to v$, hence $tv_n \to tv$ for all $t\in \R$.  Let $\exp \colon \mathcal{G} \to G$ be the exponential map.  
Then $\exp (tv_n) \to \exp (tv)$ for all $t\in \R$.  Since $H_n$ converges to $H$, $\exp (tv) \in H$ for all $t\in \R$.  This implies 
that $v\in \mathcal{H}$.  Thus, $V \subset \mathcal{H}$.

If $\dim (H_n) = \dim (H)$ for all $n$, then $\dim (\mathcal{H} _n) = \dim (\mathcal{H}) $ for all $n$. By Lemma \ref{irsl2}, 
$\dim (V) = \dim (\mathcal{H})$.  Since $V\subset \mathcal{H}$, we get that $V = \mathcal{H}$.   Thus, $\mathcal{H} _n \to \mathcal{H}$ in Sub$_\mathcal{G}$.  

\end{proof}

\section{\bf On the support of random subgroups invariant under unipotent transformations}

In this section we consider invariant measures on Grassmannian for projective linear actions and first prove a result on the support of the invariant measure.  
The result of the following kind is known as Furstenberg's Lemma (cf. Section 5, \cite{GL}) and it has been proved for 
minimally almost periodic groups in Lemma 3 of \cite{Fu} and when $\Gamma$ (as in Lemma \ref{irsl1}) is certain semisimple analytic group in Proposition 5.1 of \cite{GL}.   Here we provide the following using a result of \cite{Da}.

\begin{lemma}\label{irsl1}
Let $V$ be a finite dimensional vector space and $\Gamma$ be a subgroup of $GL(V)$. Also, let $U$ be a subgroup of $\Gamma$ consisting of unipotent transformations.  
Let $\Gr _k(V)$ be the Grassmannian manifold of $k$-dimensional subspaces of $V$ and $\lambda$ be a probability measure on 
$\Gr _k(V)$ which is invariant under the induced action of $\Gamma$ on $\Gr _k(V)$. Also, let $I_\lambda = \{ \tau \in GL(V) \mid  \tau ~~{\rm is ~~ trivial ~~ on ~} S \}$ where $S$ is the support of $\lambda$. Then $\Gamma$ is contained in an algebraic group $\tilde \Gamma$ such that $\tilde \Gamma / I_\lambda $ is a compact group, and $S$ consists of $U$-invariant subspaces. 

\end{lemma}

\begin{proof}
    
Let $\lambda $ be the $\Gamma$-invariant probability measure on $\Gr _k(V)$. 
Let $G_\lambda = \{\tau \in GL(V) \mid \tau (\lambda ) = \lambda \}$ and 
$I_\lambda = \{ \tau \in GL(V) \mid  \tau ~~{\rm is ~~ trivial ~~ on ~} S \}$ where $S$ is the support of $\lambda$. 
Then $I_\lambda $ is a normal algebraic subgroup of $G_\lambda $.  Since $\Gr _k(V)\simeq GL(V)/H$ where $H$ is the stabilizer of 
a $k$-dimensional subspace, by Corollary 2.6 of \cite{Da} we get that $G_\lambda$ is a real algebraic 
group and $G_\lambda /I_\lambda $ is compact.  Thus, any unipotent transformation in $G_\lambda $ is in $I_ \lambda$.  
This proves the result.  
\end{proof}

The following proposition is crucial to the proof of our results, it describes the action of the unipotent group on orbits that support invariant measures: for a connected Lie group $G$, Aut$(G)$ has a natural action on the Lie algebra of $G$ by identifying each $\tau \in {\rm Aut}(G)$ with its differential ${\rm d}\tau$.  

\begin{proposition}\label{irsp1}
Let $G$ be a connected Lie group with Lie algebra ${\mathcal G}$, and $S$ be a closed connected subgroup of $G$ with Lie algebra $\mathcal S$.  Let $\Gamma$ be a $\sigma$-compact subgroup of Aut$(G)$, and $U$ be a subgroup of $\Gamma$ consisting of unipotent automorphisms.  
If either 

(1) there is a $\Ga$-invariant probability measure $\lambda$ on ${\rm Sub}_{\mathcal G}$ with
$\lambda (\{ {\rm d}\tau ({\mathcal S}) \mid \tau \in \Gamma \}) =1$ 

\noindent or

(2) there is a $\Ga$-invariant probability measure $\mu$ on ${\rm Sub}_G$ with $\mu (\{ \tau (S) \mid \tau \in \Gamma \}) =1$;

then $S$ is $N$-invariant where $N$ is the smallest normal subgroup of $\Gamma$ containing $U$

\end{proposition}

\begin{proof}

Let $X= \Gr _k (\mathcal{G})$ be the Grassmannian manifold of $k$-dimensional subspaces of $\mathcal{G}$ 
where $k$ is the dimension of $\mathcal{S}$.   Consider the canonical action of $\Gamma$ on $X$ given by: 
$(\tau , V) \mapsto {\rm d}\tau (V)$ for $\tau \in \Gamma$ and $V\in X$.  

If (1) holds, then in view of Lemma \ref{irsl2}, we 
get that $\lambda$ is a $\Gamma$-invariant probability measure on $X$.  By Lemma \ref{irsl1}, $\lambda$ is 
supported on $U$-invariant subspaces.  This proves that $\mathcal{S}$ is $N$-invariant.  
Since $S$ is a connected Lie subgroup, $S$ is $N$-nvariant.

Suppose (2) holds.  Let $\mathcal {O} = \{ \tau (S) \mid \tau \in {\rm Aut}~(G ) \}$ and 
$\phi \colon  \mathcal{O} \to X$ be defined by $\phi (\tau (S) )= {\rm d}\tau (\mathcal{S})$.  Then it follows from Proposition \ref{cs} and Lemma \ref{irsl2} that $\phi$ is continuous.  Let $\lambda = \phi (\mu )$.  Then $\lambda$ is a $\Gamma$-invariant probability measure and 
$\lambda (\{ \tau ({\mathcal S}) \mid \tau \in \Gamma \}) =1$.  Now the result follows from the previous case.  
\end{proof}

\section{\bf IRS on orbits}

In this section, we show the existence of invariant random subgroups supported on various orbits of the conjugation action of a connected Lie group $G$ on $\text{Sub}_G$. In other words, we discuss the proofs of Theorems \ref{irsl}, \ref{irsK}, \ref{irsB} and \ref{irsD}. We also discuss certain structural results related to these theorems. 

\begin{proof}[Proof of Theorem \ref{irsl}]
It is sufficient to prove that (2) implies (3).  
Let $L$ be a semisimple Levi subgroup of $G$ and $N$ be the nilpotent radical of $G$ as in the statement of the theorem.  Then $L$ is a connected Lie subgroup of $G$.  Suppose that there is an IRS supported on the set $L_G$ of all semisimple Levi subgroups. Since any two semisimple Levi subgroups conjugate in $G$ by an element of $N$, there is a $N$-invariant probability measure $\mu$ such that $\mu (\{ xLx^{-1} \mid x \in N \}) =1$.  Let $\mathcal{G}$ be the Lie algebra of $G$.  Since $N$ is a nilpotent normal subgroup of $G$, Ad$(x)$ is a unipotent transformation on $\mathcal{G}$ for all $x\in N$.  Therefore, by Proposition \ref{irsp1}, $\mu$ is supported on $N$-invariant subspaces of $\mathcal{G}$.  Thus, $L$ is normalized by $N$.  Since any two semisimple Levi subgroups are conjugate by an element of $N$, we get that $L$ is normal in $G$.  

Suppose the set of all semisimple Levi subgroups is closed in $\text{Sub}_G$. Then $L_G$ is compact; and, since $N$ is amenable, 
$L_G$ admits an $N$-invariant probability measure.  Therefore, $(2)$ implies $(3)$ proves that $L$ is normal in $G$.
    
\end{proof}

\begin{proof}[Proof of Theorem \ref{irsK}]
Suppose that there is an IRS $\mu$ such that $\mu (K_G)=1$, where $K_G$ is the set of all maximal compact subgroups of $G$ with $K$ being a maximal compact subgroup of $G$. Then $K$ is connected and $K_G=\{xKx^{-1} \mid x\in G \}$ 
(cf. Theorem 3.1, Chapter XV of \cite{Ho}), hence by Proposition \ref{irsp1}, 
$K$ is normalized by any Ad-unipotent subgroup of $G$.   

Let $L$ be a semisimple Levi subgroup of $G$ and $K_l$ be a maximal compact subgroup of $L$. Then $K_l \subset gKg^{-1}$ for some $g\in G$.  Replacing $K$ by a conjugate, we may assume that $K$ contains 
$K_l$. This, in particular, implies that $K$ is normalized by compact simple factors of $L$. Let $L_1$ be a noncompact simple factor of $L$.  Then $L_1$ is generated by the Ad-unipotent elements of $L_1$. Therefore, $L_1$ normalizes $K$. Thus, the seimsimple Levi subgroup $L$ normalizes $K$.  

Now, let $N$ be the nilpotent radical of $G$. Then Ad$(x)$ is unipotent for any $x\in N$, hence $N$ normalizes $K$. Thus, $NK$ is a closed connected subgroup of $G$. Since $K$ is a maximal compact subgroup of $G$, $K$ is a maximal compact subgroup of $NK$. Since $N$ normalizes $K$, $K$ is normal in $NK$ and hence contains any compact subgroup of $NK$. Recall that any two maximal compact subgroups of a connected Lie groups are conjugate.     
Let $T=N\cap K$.  Then $T$ is the largest compact normal subgroup of $N$ and hence $N/T$ is a simply connected nilpotent Lie group.  Since $N$ is normal in $G$ and $T$ is the largest compact normal subgroup of $N$, 
$T$ is normal in $G$.  Since $T\subset K$, replacing $G$ by $G/T$, we may assume that $N$ is simply connected.  

Let $\tilde G = G/N$ and $\tilde K$ be a maximal compact subgroup of $\tilde G$ containing $KN/N$.  Then $\tilde K$ is connected and there is a closed 
subgroup $H$ of $G$ containing $N$ such that $H/N = \tilde K$.  Since $N$ and $\tilde K$ are connected, $H$ is also connected (cf. Theorem 7.4 of \cite{HeR}). As $N$ is simply connected, $H$ contains a compact subgroup $M$ such that $H= MN$ and $\tilde K= H/N \simeq M$ (cf. Theorem 3.2, Chapter XII of \cite{Ho}). Also, $K\subset H$ as $H/N =\tilde K$ contains $KN/N$. Since $M$ is a compact subgroup of $H$ and $K$ is a maximal compact subgroup, $M\subset xKx^{-1}$ for some $x\in H=MN$. But $N$ normalizes $K$, hence $M\subset K$. This shows that $\tilde K = KN/N$. 

Now, let $R$ be the solvable radical of $G$.  Then $N\subset R$ and by Corollary 3.16.4 of \cite{Va}, $R/N$ is the center of $\tilde G$, hence $\tilde K$ is normalized by $R/N$.  This shows that $KN$ is normalized by $R$.  Since the semisimple Levi subgroup $L$ normalizes $K$, $L$ normalizes $KN$ as well. It follows that $KN$ is a normal subgroup of $G$.  

As $N$ normalizes $K$, $K$ is a normal subgroup of $KN$.  
Since $N$ is a simply connected nilpotent Lie group, $N\cap K= \{  e\}$ and hence $KN/K~(\simeq N)$ has no compact subgroup.  
This proves that $K$ is a normal subgroup of $G$.  
\end{proof}
The existence of an IRS supported on the set of maximal compact subgroups has certain implications on the structure of the connected Lie group $G$ under consideration as described in the following proposition.
\begin{proposition}\label{s1}
Let $G$ be a connected Lie group with solvable radical $R$ and $M$ be a maximal compact subgroup of $G$.  Suppose $M$ is normal in $G$.  
Then,
\begin{enumerate}

\item the solvable radical $R$ admits a compact central subgroup $T$ so that $R/T$ is simply connected; 

\item any closed simple Lie subgroup of $G$ is either compact or has no nontrivial compact connected subgroups.

\end{enumerate}
\end{proposition}

\begin{remark}
The universal cover of $SL(2, \R)$ is a simple Lie group without any compact subgroup of positive dimension. Thus the second possibility in 2 of Proposition \ref{s1} may occur.
\end{remark}

\begin{proof}
Suppose the maximal compact subgroup $M$ of $G$ is normal in $G$.  If $R$ is the solvable radical of $G$, then $R\cap M$ is a 
maximal compact subgroup of $R$. Since $M$ is normal in $G$, $R\cap M$ is normal in $R$. If $T= R\cap M$, then $T$ is abelian and $R/T$ has no compact subgroup. Thus, $R/T$ is simply connected (cf. exercise 3, Chapter XII of \cite{Ho}).   

Let $S$ be any simple Lie subgroup of $G$. Suppose that there exists a non-trivial compact connected subgroup of $S$. Then $M\cap S$ is a compact normal subgroup of positive dimension in $S$.  Since $S$ is simple, $S= S\cap M$.
\end{proof}

We now apply the above results to show that a finite covolume subgroup is cocompact if the subgroup is a Levi subgroup or a maximal compact subgroup.    

\begin{corollary}\label{ak}
Let $G$ be a connected Lie group and $H$ be a closed subgroup of $G$ such that $G/H$ carries a $G$-invariant probability measure.  
Suppose $H$ is either a Levi subgroup or a maximal compact subgroup.  
Then $G/H$ is a compact group. 
\end{corollary}

\begin{proof}
We prove the result when $H$ is a Levi subgroup, the other case is similar. Suppose $H$ is a Levi subgroup of $G$. Then the stabilizer, say $G_{xH}$, at any $xH\in G/H$ is $xHx^{-1}$.  Therefore, using the map $xH\mapsto G_{xH}$, we obtain a IRS supported on the set of semisimple Levi subgroups.  
By Theorem  \ref{irsl}, we get that $G$ has a unique semisimple Levi subgroup and is normal in $G$.  This implies that $G/H$ is a locally compact group.  
Since $G/H$ carries a $G$-invariant probability measure, $G/H$ is compact.  
\end{proof}

Now, we assume that $G$ is a real algebraic group. Also, let $B$ be a Borel subgroup of $G$ and $B_G$ be the set of all Borel subgroups of $G$, i.e., $B_G=\{xBx^{-1}|x\in G\}$.

\begin{proof}[Proof of Theorem \ref{irsB}]
Since $B$ is a Borel subgroup of $G$, $G/B$ is compact, and hence there is a compact set $C$ in $G$ such that $G=CB$.  
Since any two Borel subgroups are conjugate in $G$, $B_G = \{xBx^{-1} \mid x \in C \}$.  Thus, $B_G$ is a compact subset of 
$\text{Sub}_G$.

Suppose that there is an IRS $\mu$ such that $\mu (B_G)=1$. To prove that $G$ is amenable, it is sufficient to prove that any
semisimple Levi subgroup of $G$ is compact. Let $R$ be the solvable radical of $G$. Then $R\subset B$, hence any Borel subgroup of $G$ is a lift of a Borel subgroup of $G/R$. Thus, replacing $G$ by $G/R$, we may assume that $G$ is semisimple. Now, let $S$ be a simple factor of $G$. If $S$ is non-compact, then $S$ is generated by unipotent elements. This implies by Proposition \ref{irsp1} that $B$ is normalized by $S$, hence $B\cap S$ is a solvable normal subgroup of $S$ containing a Borel subgroup of $S$.  This is a contradiction to the simplicity of $S$.  
Hence, any simple factor of $G$ is compact. This completes the proof.  
\end{proof}

\begin{proof}[Proof of Theorem \ref{irsD}]
Let $D$ be a maximal diagnolizable subgroup of $G$ and $U$ be the unipotent radical of $G$.  
Then there is a reductive Levi subgroup $L$ of $G$ containing $D$ such that $G= LU$.   
Let $L'= [L, L]$.  Then $L'$ is a semisimple Levi subgroup of $G$.   

Suppose that there is an IRS $\mu$ such that $\mu (D_G)=1$. Then, by Proposition \ref{irsp1}, 
Any unipotent subgroup of $G$ normalizes $D$. If $S$ is a noncompact simple algebraic subgroup of $G$, then the maximal split diagonalizable subgroup $D_S$ of $S$ is contained in 
a conjugate of $D$. Replacing $S$ by a conjugate of $S$, we may assume that $D_S\subset D$. Since $S$ is noncompact, $S$  is generated by unipotent elements, hence $S$ normalizes $D$. Therefore, $D\cap S$ is a normal subgroup of positive dimension which is a contradiction to the simplicity of $S$. Thus, the semisimple 
Levi subgroup $L'= [L, L]$ of $G$ is compact. This implies that $D \subset Z(L)$ because $L=[L, L] Z(L)$.  

Since the unipotent subgroups normalize $D$, $D$ is normalized by $U$.  Since $D\cap U$ is trivial, $D$ is centralized by $U$ also.  Since $G= LU$, $D$ is central in $G$.  
This proves that (1) implies (2) and that (2) implies (3) is trivial.  If $D$ is normal, the $\delta _D$, the dirac measure at $D$ is the required IRS.  This proves that (3) implies (1).    

Suppose (1) holds.  Then the proof of the first implication shows that Levi subgroup is compact and $D$ is central.  Let $R$ be the solvable radical of $G$. Then $D_R=R\cap D$ is a central subgroup of $R$, hence $D_RU$ is a nilpotent subgroup of $R$ where $U$ is the unipotent radical of $G$.  Since $R/D_RU$ is an anisotropic torus (that is, without any split torus subgroup) we get that $R/D_RU$ is compact. 
\end{proof}  

\begin{corollary}\label{a2}
Let $G$ be a connected real algebraic group and $D$ be a maximal diagonalizable subgroup of $G$. Suppose $G/D$ carries a $G$-invariant probability measure. Then $G$ is a direct product of $D$ and a compact group. 
\end{corollary}

\begin{proof}
Let $m$ be the $G$-invariant probability measure on $G/D$.  Then stabilizer, $G_{xD}$, 
at any $xD\in G/D$ is $xDx^{-1}$.  Therefore, the map $xD \mapsto xDx^{-1}$ gives rise to an IRS supported on maximal diagonalizable subgroups.  
By Theorem \ref{irsD}, we get that $D$ is normal in $G$ and the Levi subgroup of $G$ is compact.    
Since $G/D$ has a $G$-invariant probability measure, $G/D$ is a compact group.  Therefore, the unipotent radical of $G$ is trivial.  Thus, 
$G$ is a direct product of $D$ and a compact subgroup.   
\end{proof}
\acknowledgement{ Both the authors would like to thank the anonymous referee who suggested the splitting of the previous version of the article, and many other valuable suggestions on improving the exposure of the article especially suggesting \cite{A7s0}. The authors also like to thank International Centre for Theoretical Sciences, Bangalore for its hospitality during the programme Ergodic Theory and Dynamical Systems, $5$th December to $16$th December, 2022 when this work was initiated.
The first author thanks the Indian Statistical Institute Bangalore Centre for its hospitality during his several visits while this work was going on. The second author acknowledges the support from SERB through the grant under MATRICS MTR/2022/000429.}

\bibliographystyle{plain}
\bibliography{refnil}

\end{document}